\documentclass[a4,semcolor,psfig,english,11pt,epsf,portrait]{article}

\usepackage{amsmath,amsfonts,latexsym,amscd,amssymb,theorem}
\usepackage{bigints,moreverb,rotating,graphics,amsfonts,latexsym,amscd}
\usepackage[ansinew]{inputenc}
\usepackage{epsfig}
\usepackage{relsize}
\usepackage{url}
\usepackage{psfrag}
\usepackage[usenames,dvipsnames]{color}
\usepackage{color}
\usepackage{mathdots}
\usepackage{hyperref}
\hypersetup{colorlinks=true,allcolors=blue}
\usepackage[T1]{fontenc}
\usepackage{epsfig}
\usepackage{graphicx}
\usepackage[T1]{fontenc}

\newcommand{\BE}{\begin{equation}}
\newcommand{\EE}{\end{equation}}
\newcommand{\Bii}{\begin{eqnarray*} }
\newcommand{\Eii}{\end{eqnarray*} }
\newcommand{\BEii}{\begin{eqnarray} }
\newcommand{\EEii}{\end{eqnarray} }
\newcommand{\BA}{\begin{array} }
\newcommand{\EA}{\end{array} }
\def\R{\mathbb{R}}

\def\<{\langle}
\def\>{\rangle}

\def\p{\partial}
\def \diam{\hbox{\rm diam }}

\newtheorem{thm}{Theorem}[section]

\newtheorem{theorem}[thm]{Theorem}

\newtheorem{proposition}[thm]{Proposition}

\newtheorem{rem}[thm]{Remark}
\newenvironment{preuve}
{\medskip\noindent{\bf Proof.\\ }}{\null \hfill
$\blacksquare$\par\medskip}
\newcommand{\x}{{\bf x}}

\newcommand*\colvec[3][]{
    \begin{pmatrix}\ifx\relax#1\relax\else#1\\\fi#2\\#3\end{pmatrix}
}

\def\R{\mathbb{R}}
\def\N{\mathbb{N}}

\def\<{\langle}
\def\>{\rangle}

\begin{document}
\title{Inverse moving point source problem for the wave equation}

\author{
\normalsize\textsc{Hanin Al Jebawy $^{1}$, Abdellatif El Badia$^{1}$ and  Faouzi Triki$^2$}} \vspace{10pt}
\maketitle \footnotetext[1] {\small Sorbonne universit\'es, Universit\'e de Technologie de Compi\`egne, Laboratoire de Math\'ematiuqes Appliqu\'ees de
Compi\`egne, Centre de recherche Royallieu,
CS 60 319, 
60203 Compi\`egne cedex, France. Emails:  hanin.aljebawy@utc.fr,\,  abdellatif.elbadia@utc.fr}
\maketitle 

\footnotetext[2]{\small Universit\'e Grenoble Alpes, Laboratoire Jean Kuntzmann -B\^{a}timent IMAG, 700 Avenue Centrale 38401 Saint-Martin-d'H\`{e}res, France. Emails: faouzi.triki@univ-grenoble-alpes.fr\\The authors were supported by the grant ANR-17-CE40-0029 of the French National Research Agency ANR (project MultiOnde)}
\begin{abstract}
 In this paper, we consider the problem of identifying a single moving point  source for a three-dimensional wave equation from boundary measurements. Precisely,  we show that the knowledge of the field  generated by the source at six different points of the boundary over a finite time interval is sufficient to determine  uniquely its  trajectory. We also derive a Lipschitz stability estimate for the inversion. 
\end{abstract}

\section{Introduction}
 
Inverse source problems are of importance  in several  scientific areas including biomedical engineering,  antenna synthesis, geology, and medical imaging   \cite{Ammari02, Bao10, Bao11, Bao21, Bao212, Badia2002, Hu, Isakov90,Liu}.  In this paper we consider the inverse source problem for the  wave equation. Precisely, we study the problem of determining  the trajectory of  a moving source in a bounded domain from a single boundary measurement.  Identification of sources with time-varying  locations  has many significant applications such as the recovery  of  mobile pollution sources, or small debris in low-earth orbit, and underwater sonar systems.
\medskip\\
 We assume here  that  media surrounding the  point source  is homogeneous and isotropic, and the measurement of the wave field is provided only  on a  small part of the  boundary of the domain. Let $\phi$ to be the field generated by a single point source, that is, a solution of the following initial value problem for the three dimensional wave equation
\begin{equation}
\left\{\begin{aligned}\label{phi}
&\frac{1}{c^2}\phi_{tt}-\Delta\phi=\lambda\delta({\bf x}-b(t))&\text{in }&\R^3\times (0,T)\\
& \phi({\bf x},0)=\phi_t({\bf x},0)=0&\text{in }&\R^3,
\end{aligned}
\right.
\end{equation}
where $T>0$ is a fixed time,  $c>0$ is the speed of the wave, $\lambda>0$ is the intensity, and   $b\in C^2([0,T]; \mathbb R^3)$ is the position of the point source  confined within a bounded  domain $D\subset\R^3$. Let $\Omega$ be a smooth bounded domain satisfying  $\overline D\subset \Omega$ with boundary $\Gamma$. Notice that the trajectory of the point source  remains away from $\Gamma$.  Define  $G_b(T)$ to be the graph of the function $b$, that is,   $G_b(T)=  \{(s, b(s)); \, s\in [0, T] \} $.   We also assume that  source is subsonic, in other words,  the speed of the 
source which is  the first derivative of $b$, satisfies  
 \BE \label{subsonic}
  \|b^\prime\|_{C([0,T])}\leq  c_0<c.
 \EE

The main goal of this paper is to reconstruct the trajectory followed by the source  $b$ by measuring  $\phi$ on a part of the surface $\Gamma$.
 \medskip\\
 
There are few  works dealing with  inverse moving source problems for  the wave equation, and different approaches have been used for solving it. In \cite{Badia2001,Ohe2011}, for example, the authors considered the problem with time-varying point sources, and they applied some algebraic direct methods for reconstructing these sources from a single boundary measurement, based on the concept of the reciprocity gap functional. Later, in \cite{Ohe2020} the author provided a generalization of this algorithm to the problem of moving sources, yet this method includes some difficulties concerning the assumptions on the trajectories followed by the sources, and no stability estimate has been provided. Another algebraic algorithm for the reconstruction of one moving source was given by \cite{Nakaguchi2012} using the observed values of the retarded potential and all its derivatives at a single observation point. Other optimization techniques were also used in solving this problem, see for example \cite{Bruncker2000,Komornik2002,Komornik2005,Rashedi2015}.
 \medskip\\
 In our work, we are interested  in identifying the trajectory $b(t)$ of a single moving point source source $\delta(\cdot-b(t))$ with a known intensity $\lambda>0$, from the measurement of the generated field $\phi$ at six well-chosen points $\x_i, \, i=1,\cdots, 6,$ located on  the observation surface $\Gamma$.
 \medskip\\  
 Throughout the paper we denote   $\cdot$, and $|\cdot|$  the  scalar product and the Euclidean norm respectively in
 $\mathbb R^3$.  
  \medskip\\  
The unique  solution $\phi \in C([0, T]; L^2_{\textrm{Loc}}(\mathbb R^3))\cap C^1([0, T]\times (\mathbb R^3\setminus D))$ of \eqref{phi} is given by \cite{Jackson99,Nakaguchi2012,Ohe2011}
\begin{equation}\label{sln}
\phi({\bf x},t)=\frac{\lambda}{4\pi c}\frac{ Y(r)}{\ |{\bf x}-b(r)|\ h({\bf x},r)},
\end{equation}
where $Y$ denotes the Heaviside function, $r \in C([0, T]\times \overline \Omega) \cap C^1\left(([0, T]\times \overline \Omega)\setminus G_b(T) \right) $ is the unique solution to the equation 
\begin{equation}\label{r}
r({\bf x},t)=t-\frac{|{\bf x}-b(r)|}{c},
\end{equation}
for each fixed $(\x, t)$, and 
$$h({\bf x},r)=1-\frac{b^\prime(r).({\bf x}-b(r))}{c|{\bf x}-b(r)|}.$$
We note that since $b\in C^2([0,T], \mathbb R^3)$ satisfies \eqref{subsonic}, we have 
\begin{equation}\label{c0}
	h\geq1- c^{-1}\|b^\prime\|_{C([0,T])} :=h_0>0.
\end{equation}
 Moreover, differentiating \eqref{r} in $t$, we deduce that 
 \begin{equation}\label{odeR}
 \frac{\p r}{\p t}({\bf x},t)=\frac{1}{h({\bf x},r)}>0,
 \end{equation}
  and so  $r$ is strictly increasing in $t$.
\medskip\\ 
The objective of our work is to prove that a single observation $\phi$ on the surface $\Gamma$ uniquely determines the source term $b$. Our strategy is to first prove that the boundary observation $\phi$ uniquely determines $r$, which in return allows us to reconstruct $b$ using the relation \eqref{r}.
Our goal thus is to reconstruct $r$ for different positions ${\bf x}$ on $\Gamma$. First, we assume that
\begin{equation}\label{T}
T>\underset{{\bf x}\in \Gamma, {\bf y}\in D}{\sup}\frac{|{\bf x}-{\bf y}|}{c}:= T_0.
\end{equation}
Since the wave field is propagating with a finite speed $c$, the assumption \eqref{T} is sufficient to allow the 
information on the source to arrive to the  observation surface $\Gamma$.\\

Now, we take ${\bf x}$ on $\Gamma$, fixed but arbitrary, and we define the time 
\begin{equation}\label{tx}
t_{{\bf x}}=\sup\{t>0;\ \phi({\bf x},t)=0\}.
\end{equation}
We note that  from \eqref{sln}, and the definition of the Heaviside function, we have $\phi>0$ for $r>0$. Due to 
the initial conditions the set $\{t>0;\ \phi({\bf x},t)=0\}$ is not empty,  and  considering \eqref{r} we remark that $r({\bf x},t)>0$ holds  for $t$ large enough. Thus, for $T$ satisfying \eqref{T}, we guarantee that $\phi({\bf x}, .)$ is not zero on $(0,T)$, and thus $t_{\bf x}\in (0,T)$. \\\

 On the other hand since $r$ is continuous and strictly increasing in $t$, we deduce from \eqref{sln} that $r({\bf x},t_{{\bf x}})=0$. Therefore, using the relation \eqref{r}, we get
 \begin{equation} \label{relationtx}
 |{\bf x}-b(0)|=ct_{{\bf x}}.
 \end{equation}
Repeating the same procedure for different ${\bf x}\in \Gamma$, one can estimate the location of $b(0)$.
\medskip\\
Our goal now is to reconstruct $r$ for different ${\bf x}$ on $\Gamma$, which allows us later to reconstruct $b$ using the relation \eqref{r}.
\medskip\\
Our paper is organized as follows: First, we provide in section \ref{sec2} an ODE based method for the reconstruction of $r$ on $\Gamma$, then we reconstruct in section \ref{sec3} the trajectory followed by the source $\{b(t); t\in [0, T]\}$ using the previously calculated values of $r$ and the measurements of $\phi$ on six well chosen observation points on $\Gamma$. The uniqueness of the reconstruction is announced in Theorem \ref{uniqueness}.
 Finally, stability estimates are derived  in section \ref{sec4}. The stability in the recovery of the trajectory is provided 
 in Theorem \ref{stabilityestimate}.
 
\section{Reconstruction of $r$}\label{sec2}
Let ${\bf x}\in \Gamma$, fixed but arbitrary. Since $r$ satisfies \eqref{odeR}, we deduce from \eqref{sln} that for every $t>t_{{\bf x}}$, defined in \eqref{tx}, $r$ satisfies the following equation
\begin{equation}\label{ODE}
\left\{\begin{aligned}
&\frac{\p r}{\p t}({\bf x},t)+\frac{4\pi c \phi({\bf x},t)}{\lambda}r= \frac{4\pi c \phi({\bf x},t)}{\lambda} t\\
& r({\bf x},t_{{\bf x}})=0.
\end{aligned}
\right.
\end{equation}
Since $\phi(\x, t)$ is given on the boundary our goal here  is to solve \eqref{ODE} in $r$ for different ${\bf x}$ on $\Gamma$.
\begin{rem}We note that $r$ is a strictly increasing function in $t$. Thus, as $t$ varies between $t_{{\bf x}}$ and $T$ we have $0\leq r({\bf x},t)\leq r({\bf x},T)$, where $r({\bf x},T)=T-\frac{|{\bf x}-b(r({\bf x},T))|}{c}$. Therefore, in order to guarantee the reconstruction of $b(t)$ for $t\in (0,T)$, we  may assume that our observations continue until a later time 
\begin{equation}\label{T1}
T_{\textrm{obs}}=T+T_0\geq  T+ \underset{{\bf x}\in \Gamma}{\sup}\frac{|{\bf x}-b(T)|}{c},
\end{equation}
where $T_0$ is defined in \eqref{T}. 
\end{rem}
For ${\bf x}\in \Gamma$ fixed, we consider the system \eqref{ODE} with $t\in (t_{\bf x},T_{\bf x})$ where $T_{\bf x}= T+\frac{|{\bf x}-b(T)|}{c}$. 

\begin{proposition}
The system \eqref{ODE} has a unique solution $r\in C^1([t_\x, T_\x]; \mathbb R_+)$.  
\end{proposition}
 \begin{preuve}
Since $\phi\in  C^1([0, T_{\textrm{obs}}]\times (\mathbb R^3\setminus D))$, the result is a direct consequence  of Cauchy-Lipschitz Theorem. 
\end{preuve}

 Now, we are able to reconstruct $|{\bf x}-b(r({\bf x},t))|$ for every $t\in (t_{{\bf x}}, T_{\bf x})$. Since $r$ is smooth and strictly increasing in $t$, thus we can reconstruct $|{\bf x}-b(r)|$ for every $r\in (0,T)$. Repeating the same procedure for finite number of  positions ${\bf x}_i$ on $\Gamma$, one can estimate the location of $b(r)$. This will be explained in details in the next section.
\section{Reconstruction of $b$}\label{sec3}

The goal of this section is to reconstruct the trajectory followed by the source $b$ using the previously calculated values of $r$ on the observation surface $\Gamma$. Our algorithm requires the measurement of $\phi$ on well chosen observation points $\{{\bf x_i}\}_{1\leq i\leq 6}$ on $\Gamma$.
\medskip\\
First, we define for every ${\bf x_i}\in \Gamma$
$$t_{\bf x_i}= \frac{|{\bf x_i}-b(0)|}{c} \quad \text{and} \quad T_{\bf x_i}=T+\frac{|{\bf x_i}-b(T)|}{c}.$$
 We assume in this section that the functions $r({\bf x_i},t)$ for $t\in (t_{\bf x_i}, T_{\bf x_i})$ are previously constructed. Now we take $\tau\in (0,T)$, fixed but arbitrary, and ${\bf x_1}=(x_1,y_1,z_1)\in \Gamma$, then from the regularity of $r$ we deduce that there exists $t_{1,\tau}\in (t_{\bf x_1},T_{\bf x_1})$ such that
$$r({\bf x_1},t_{1,\tau})=\tau.$$
 Thus, we deduce from \eqref{r} that $b(\tau)$ satisfies
$$|{\bf x_1}-b(\tau)|=c(t_{1,\tau}-\tau).$$
Therefore, $b(\tau)$ moves on a sphere of center ${\bf x_1}$ and radius $c(t_{1,\tau}-\tau)$, which implies that $b(\tau)=(b_1(\tau),b_2(\tau),b_3(\tau))$ satisfies the equation
\begin{equation}\label{x1}
(b_1(\tau)-x_1)^2+(b_2(\tau)-y_1)^2+(b_3(\tau)-z_1)^2=c^2(t_{1,\tau}-\tau)^2.
\end{equation}
Furthermore, taking ${\bf x_2}=(x_2,y_2,z_2)$ another point on $\Gamma$, and repeating the previous procedure, we know that $b(\tau)$ moves on a sphere of center ${\bf x_2}$ and radius $c(t_{2,\tau}-\tau)$ for some $t_{2,\tau}\in (t_{\bf x_2},T_{\bf x_2})$ that satisfies
$$r({\bf x_2},t_{2,\tau})=\tau.$$
Therefore, $b(\tau)$ also satisfies
\begin{equation}\label{x_2}
(b_1(\tau)-x_2)^2+(b_2(\tau)-y_2)^2+(b_3(\tau)-z_2)^2=c^2(t_{2,\tau}-\tau)^2.
\end{equation}
Subtracting \eqref{x1} and \eqref{x_2} we get
\begin{equation}
2(x_2-x_1)b_1(\tau)+2(y_2-y_1)b_2(\tau)+2(z_2-z_1)b_3(\tau)\hspace{-1mm}=
\hspace{-1mm}|{\bf x_2}|^2-|{\bf x_1}|^2\hspace{-1mm}+c^2(t_{1,\tau}-\tau)^2\hspace{-1mm}-c^2(t_{2,\tau}-\tau)^2.\label{12}
\end{equation}
Similarly, taking ${\bf x_3}, {\bf x_4},{\bf x_5}$ and ${\bf x_6}\in\Gamma$, we get
 \begin{equation}\label{34}
 2(x_4-x_3)b_1(\tau)+2(y_4-y_3)b_2(\tau)+2(z_4-z_3)b_3(\tau)\hspace{-1mm}=\hspace{-1mm}|{\bf x_4}|^2\hspace{-1mm}-|{\bf x_3}|^2\hspace{-1mm}+c^2(t_{3,\tau}-\tau)^2\hspace{-1mm}-c^2(t_{4,\tau}-\tau)^2,
 \end{equation}
 and
 \begin{equation}\label{56}
 2(x_6-x_5)b_1(\tau)+2(y_6-y_5)b_2(\tau)+2(z_6-z_5)b_3(\tau)\hspace{-1mm}=|{\bf x_6}|^2\hspace{-1mm}-|{\bf x_5}|^2+c^2(t_{5,\tau}-\tau)^2\hspace{-1mm}-c^2(t_{6,\tau}-\tau)^2.
 \end{equation}
 Equations \eqref{12}, \eqref{34}, and \eqref{56} can be rewritten in the matrix form
 $$XB=\frac{1}{2}A,$$
 where 
 \begin{equation}\label{X}
 X=\left(\begin{array}{ccc}
 x_2-x_1 & y_2-y_1 & z_2-z_1\\
 x_4-x_3 & y_4-y_3 & z_4-z_3\\
 x_6-x_5 & y_6-y_5 & z_6-z_5
 \end{array}\right),
 \end{equation}
 $$B=\left(\begin{array}{c}
 b_1(\tau)\\
 b_2(\tau)\\
 b_3(\tau)
 \end{array}\right),$$
 and
\begin{equation}\label{A}
 A= \left(\begin{array}{c}
 A_{12}\\
 A_{34}\\
 A_{56}
 \end{array}
 \right),
\end{equation}
 with 
\begin{equation}\label{Aij}
	A_{ij}=|{\bf x_j}|^2-|{\bf x_i}|^2+c^2(t_{i,\tau}-\tau)^2-c^2(t_{j,\tau}-\tau)^2.
\end{equation}
 Now, since $\Gamma$ is a boundary of a connected  domain one can choose the observation points $\{{\bf x_i}\}_{1\leq i\leq 6}$ such that $X$ becomes invertible, then $B$ can be reconstructed as
 \begin{equation}\label{B}
 B=\frac{1}{2}X^{-1}A.
 \end{equation}
 Notice that $X$ depends only on $\Gamma$.   Finally, we give a brief summary of the reconstruction of $b$ that passes through five mains steps.
 \medskip\\
 {\bf Step 1.} Choose $\{{\bf x_i}\}_{1\leq i\leq 6}$ on $\Gamma$ such that the matrix $X$ defined in \eqref{X} is invertible.
 \medskip\\
 {\bf Step 2.} Calculate the values of $t_{\bf x_i}$ given by
 $$t_{\bf x_i}=\sup\{t>0;\ \phi({\bf x_i},t)=0\}.$$
 {\bf Step 3.} Construct $r({\bf x_i},t)$ as solutions of the  equations \eqref{ODE}.
 \medskip\\
 {\bf Step 4.} For every $\tau\in (0,T)$, evaluate $t_{i,\tau}$ given by
 $$r({\bf x_i},t_{i,\tau})=\tau.$$
 {\bf Step 5.} Construct the matrix $A$ given by \eqref{A}, and finally evaluate $b(\tau)$ using the relation \eqref{B}.\\
 
 Then, we have the following result. 
 \begin{theorem} \label{uniqueness}
Let $\{{\bf x_i}\}_{1\leq i\leq 6}$ be fixed points  on $\Gamma$ chosen such that the matrix $X=(\x_2-\x_1, \x_4-\x_3, \x_6-\x_5)^T$ is invertible. Then the knowledge of $\phi(\x_i, t), \, i=1, \cdots, 6$ for $t\in [0, T_{\textrm{obs}}]$ where $T_{\textrm{obs}}= T+T_0$, determines uniquely the trajectory of the  point source, that is  $\{b(t); t\in[0, T]\}$.
 \end{theorem}
 \section{Stability estimates}\label{sec4}
 In this section we give a stability estimate for the reconstruction of $b$ from the measurements of $\phi$ on the observation points $\{{\bf x_i}\}_{1\leq i\leq 6}$. Our work is divided into three steps. In fact, we notice from the previous section that the reconstruction of $b$ passes through three main stages: the reconstruction of $t_{\bf x_i}$, $r({\bf x_i},t)$ and finally the trajectory $b(t)$. For this purpose, we first give the stability estimates for $t_{\bf x_i}$ and $r$, then we deduce that of the source $b$.
 \medskip\\
 Recall that $\Omega\subset\R^3$ is a bounded domain of boundary $\Gamma$, satisfying  $\overline{D}\subset \Omega$. We consider as in the previous section a set of observation points $\{{\bf x_i}\}_{1\leq i\leq 6}$ such that the matrix $X$ defined in \eqref{X} is invertible, and we suppose that we have two observations $\phi$ and $\widetilde{\phi}$ corresponding respectively to two trajectories $b$ and $\tilde b$.
 Then, we establish the following stability estimates.
 \begin{theorem}[Stability estimate of $t_{\bf x_i}$]\label{testimate}
 Let $\{{\bf x_i}\}_{1\leq i\leq 6}\in \Gamma$ such that the matrix \eqref{X} is invertible, and take $t_{\bf x_i}$ and $\tilde{t}_{\bf x_i}$ defined in \eqref{tx} corresponding to the observations $\phi$ and $\widetilde{\phi}$ respectively. Then, the following estimate holds.
 \begin{equation}\label{t bound}
 |t_{\bf x_i}-\tilde{t}_{\bf x_i}|\leq \frac{8\pi T_{\textrm{obs}} \diam(\Omega)}{\lambda}\|\widetilde{\phi}({\bf x_i},.)-\phi({\bf x_i},.)\|_{L^{\infty}(0,T_{\textrm{obs}})}, \quad i=1, \cdots, 6,
 \end{equation}
 where $T_{\textrm{obs}}= T+T_0$, with $T_0$ defined in \eqref{T}.
 \end{theorem} 
 \begin{preuve}
 We assume without loss of generality that $\tilde{t}_{\bf x_i}<t_{\bf x_i}$, then we have
 \begin{align*}
 \int_{\tilde{t}_{\bf x_i}}^{t_{\bf x_i}}\widetilde{\phi}({\bf x_i},s)ds&=\int_{\tilde{t}_{\bf x_i}}^{t_{\bf x_i}}\widetilde{\phi}({\bf x_i},s)-\phi({\bf x_i},s)ds\\
 &\leq \int_{\tilde{t}_{\bf x_i}}^{t_{\bf x_i}}|\widetilde{\phi}({\bf x_i},s)-\phi({\bf x_i},s)|ds\\
 &\leq \int_{0}^{T_{\textrm{obs}}}|\widetilde{\phi}({\bf x_i},s)-\phi({\bf x_i},s)|ds\\
 &\leq T_{\textrm{obs}} \|\widetilde{\phi}({\bf x_i},.)-\phi({\bf x_i},.)\|_{L^{\infty}(0,T_{\textrm{obs}})}.
 \end{align*}
 Therefore,
 $$\min \widetilde{\phi}({\bf x_i},.)\ (t_{\bf x_i}-\tilde{t}_{\bf x_i})\leq T_{\textrm{obs}}\|\widetilde{\phi}({\bf x_i},.)-\phi({\bf x_i},.)\|_{L^{\infty}(0,T_{\textrm{obs}})}. $$
 Furthermore, for every $t>\tilde{t}_{\bf x_i}$ we have
 $$\widetilde{\phi}=\frac{\lambda}{4\pi |{\bf x_i}-\tilde{b}(\tilde{r})|\tilde{h}({\bf x_i},\tilde{r})},$$
 where,
 $\tilde{h}({\bf x_i},\tilde{r})=1-\frac{\tilde{b}'(\tilde{r}).({\bf x_i}-\tilde{b}(\tilde{r}))}{c|{\bf x_i}-\tilde{b}(\tilde{r}|}\leq 1+\frac{\|\tilde{b}'\|_\infty}{c}$. Moreover, since $\|\tilde{b}'\|_\infty<c$, we deduce that $\tilde{h}({\bf x_i},\tilde{r})<2$. Therefore,
 $$\tilde{\phi}\geq\frac{\lambda}{8\pi\diam(\Omega)}.$$
 Finally, we obtain  
 $$(t_{\bf x_i}-\tilde{t}_{\bf x_i})\leq \frac{8\pi T_{\textrm{obs}}\diam(\Omega)}{\lambda}\|\widetilde{\phi}({\bf x_i},.)-\phi({\bf x_i},.)\|_{L^{\infty}(0,T_{\textrm{obs}})}.$$
 \end{preuve}
 \begin{theorem}[Stability estimate for $r({\bf x_i},.)$]\label{r estimate}
 Let $\{{\bf x_i}\}_{1\leq i\leq 6}\in \Gamma$ such that the matrix \eqref{X} is invertible, and take $r({\bf x_i},t)$ and $\tilde{r}({\bf x_i},t)$ solutions of \eqref{ODE} corresponding to the observations $\phi$ and $\widetilde{\phi}$ respectively with $r({\bf x_i},t_{{\bf x_i}})=0$ and $\tilde{r}({\bf x_i},\tilde{t}_{{\bf x_i}})=0$. Choose $t_{0,i}=\min\{t_{{\bf x_i}}, \tilde{t}_{{\bf x_i}}\}$, then, the following estimate holds:
 $$\|r({\bf x_i},.)-\tilde{r}({\bf x_i},.)\|_{L^{\infty}(t_{0,i},T_{\textrm{obs}})}\leq C\|\phi-\tilde{\phi}\|_{L^\infty(0,T_{\textrm{obs}})},$$
 for some $C= C(\Omega,\lambda,T, c, c_0, D)>0$.
 \end{theorem}
 \begin{preuve}
 Further $C$ denotes a  generic strictly positive constant that depends on $(\Omega,\lambda,T, c, c_0, D)$, and 
which may be different from line to line. \\
 Let ${\bf x_i}\in \Gamma$ for $1\leq i\leq 6$. Without loss of generality, we assume that $t_{{\bf x_i}}\geq \tilde{t}_{{\bf x_i}}$, our proof then is divided into two steps. First, we consider the case $t> t_{{\bf x_i}}$, then we have
 $$\p_t(r({\bf x_i},t)-\tilde{r}({\bf x_i},t))=\frac{4\pi c}{\lambda}(\phi({\bf x_i},t)-\widetilde{\phi}({\bf x_i},t))t-\frac{4\pi c}{\lambda}\phi({\bf x_i},t)r({\bf x_i},t)+\frac{4\pi c}{\lambda}\widetilde{\phi}({\bf x_i},t)\tilde{r}({\bf x_i},t). $$
 Thus,
  $$\p_t(r({\bf x_i},t)-\tilde{r}({\bf x_i},t))=\frac{4\pi c}{\lambda}(\phi-\widetilde{\phi})t-\frac{4\pi c}{\lambda}\phi(r-\tilde{r})-\frac{4\pi c}{\lambda}(\phi-\widetilde{\phi})\tilde{r}.$$
  Integrating between $t_{\bf x_i}$ and $t$, we get
  \begin{eqnarray*}
  & r({\bf x_i},t)-\tilde{r}({\bf x_i},t)=&\tilde{r}({\bf x_i},\tilde{t}_{\bf x_i})-\tilde{r}({\bf x_i},t_{\bf x_i})\\&&+\frac{4\pi c}{\lambda}\int_{t_{\bf x_i}}^{t}(\phi-\widetilde{\phi})s\ ds-\frac{4\pi c}{\lambda}\int_{t_{\bf x_i}}^{t}\phi(r-\tilde{r})ds-\frac{4\pi c}{\lambda}\int_{t_{\bf x_i}}^{t}(\phi-\widetilde{\phi})\tilde{r}ds.
 \end{eqnarray*}
  Therefore,
\begin{align*}
 |r({\bf x_i,t})-\tilde{r}({\bf x_i},t)|\leq   \|\tilde{r}'\|_{L^\infty(0,T_{\textrm{obs}})}|\tilde{t}_{\bf x_i}-t_{\bf x_i}|+\frac{4\pi c}{\lambda}T_{\textrm{obs}}\|\phi-\widetilde{\phi}\|_{L^\infty(0,T_{\textrm{obs}})}(t-t_{\bf x_i})\\
 +\frac{4\pi c}{\lambda}\|\phi\|_{L^\infty(0,T_{\textrm{obs}})}\int_{t_{\bf x_i}}^{t} r-\tilde{r}\ ds+\frac{4\pi c}{\lambda}\|\tilde{r}\|_{L^\infty(0,T_{\textrm{obs}})}\|\phi-\widetilde{\phi}\|_{L^\infty(0,T_{\textrm{obs}})}(t-t_{\bf x_i}).
\end{align*}
Using \eqref{t bound}, we deduce that
 \begin{eqnarray*}
 |r({\bf x_i,t})-\tilde{r}({\bf x_i},t)|\leq \frac{8\pi T_{\textrm{obs}}\diam(\Omega)}{\lambda h_0}\|\phi-\widetilde{\phi}\|_{L^\infty(0,T_{\textrm{obs}})}\\+\frac{8\pi c}{\lambda}T_{\textrm{obs}}^2\|\phi-\widetilde{\phi}\|_{L^\infty(0,T_{\textrm{obs}})}+\frac{c}{\textrm{dist}(\Gamma,D)h_0}\int_{t_{\bf x_i}}^{t} r-\tilde{r}\ ds,
  \end{eqnarray*}
where $h_0$ is defined in \eqref{c0}. Therefore,
$$|r({\bf x_i,t})-\tilde{r}({\bf x_i},t)|\leq C\|\phi-\widetilde{\phi}\|_{L^\infty(0,T_{\textrm{obs}})}+\frac{c}{\textrm{dist}(\Gamma,D)h_0}\int_{t_{\bf x_i}}^{t} r-\tilde{r}\ ds.$$
Applying Gronwall's Lemma we deduce that
\begin{equation}\label{case1}
\begin{aligned}
|r({\bf x_i,t})-\tilde{r}({\bf x_i},t)|&\leq C\|\phi-\widetilde{\phi}\|_{L^\infty(0,T_{\textrm{obs}})}e^{\frac{c}{\textrm{dist}(\Gamma,D)h_0}(t-t_{\bf x_i})}\\
&\leq Ce^{\frac{cT_{\textrm{obs}}}{\textrm{dist}(\Gamma,D)h_0}}\|\phi-\widetilde{\phi}\|_{L^\infty(0,T_{\textrm{obs}})}.
\end{aligned}
\end{equation}
The second case is for $t\in (\tilde{t}_{\bf x_i},t_{\bf x_i})$, then in this case we choose $t^*>t_{\bf x_i}$ such that $$|t^*-t|\leq |t_{\bf x_i}-\tilde{t}_{\bf x_i}|.$$
Then, we get
\begin{align*}
|\tilde{r}({\bf x_i,t})-r({\bf x_i},t)|&\leq |\tilde{r}({\bf x_i,t})-\tilde{r}({\bf x_i,t^*})|+|\tilde{r}({\bf x_i,t^*})-{r}({\bf x_i,t^*})|+|{r}({\bf x_i,t^*})-r({\bf x_i},t)|\\
& \leq \|\tilde{r}'\|_{L^\infty(0,T_{\textrm{obs}})}|t^*-t|+|\tilde{r}({\bf x_i,t^*})-{r}({\bf x_i,t^*})|+\|r'\|_{L^\infty(0,T_{\textrm{obs}})}|t^*-t|\\
& \leq \frac{2}{h_0}|t_{\bf x_i}-\tilde{t}_{\bf x_i}|+|\tilde{r}({\bf x_i,t^*})-{r}({\bf x_i,t^*})|. 
\end{align*}
Therefore, following \eqref{case1}, and the results of Theorem \ref{testimate}, we finally deduce that
$$|\tilde{r}({\bf x_i,t})-r({\bf x_i},t)|\leq\left(\frac{16\pi T_{\textrm{obs}}\diam(\Omega)}{\lambda h_0}+Ce^{\frac{cT_{\textrm{obs}}}{\textrm{dist}(\Gamma,D)h_0}}\right)\|\phi-\widetilde{\phi}\|_{L^\infty(0,T_{\textrm{obs}})}.$$
 \end{preuve}
\begin{theorem}[Stability estimate for b] \label{stabilityestimate}
Consider the two sources trajectories $b$ and $\tilde{b}$ with two different observations $\phi$ and $\widetilde{\phi}$ at the  points $\{{\bf x_i}\}_{1\leq i\leq 6}\in \Gamma$, such that the matrix \eqref{X} is invertible, then the following stability estimate holds 
\begin{eqnarray}
	&\|b-\tilde{b}\|_{L^\infty(0,T)}\leq &\nonumber
	 \\& C \hspace{-1mm}\underset{i=1,3,5}{\sup}\hspace{-2mm}\left(\|\widetilde{\phi}({\bf x_i,.})-\phi({\bf x_i,.})\|_{L^\infty(0,T_{\textrm{obs}})}\hspace{-1mm}+
	\hspace{-1mm}\|\widetilde{\phi}({\bf x_{i+1},.})-\phi({\bf x_{i+1},.})\|_{L^\infty(0,T_{\textrm{obs}})}\right), \nonumber\\
	\label{b estimate} 
\end{eqnarray}
 for some $C= C(\Omega,\lambda,T, c, c_0, D)>0$.
\end{theorem}
\begin{preuve}
 $C$ denotes a  generic strictly positive constant that depends on $(\Omega,\lambda,T, c, c_0, D)$, and 
which may be different from line to line. \\
Let $\tau\in(0,T)$, then following the work done in section \ref{sec3} we deduce that
\begin{equation}\label{b}
b(\tau)-\tilde{b}(\tau)=\frac{1}{2}X^{-1}(A-\widetilde{A}),
\end{equation}
where $X$ is defined in \eqref{X}, and the matrices $A$ and $\widetilde{A}$ are defined in \eqref{A} with observations $\phi$ and $\widetilde{\phi}$ respectively. Moreover, following \eqref{Aij},  we  obtain 
\begin{align*}
	A_{i\ i+1}-\widetilde{A}_{i\ i+1}=c^2\left[(t_{i,\tau}-\tau)^2-(\tilde{t}_{i,\tau}-\tau)^2\right]-c^2\left[(t_{i+1,\tau}-\tau)^2-(\tilde{t}_{i+1,\tau}-\tau)^2\right]\\
	=c^2(t_{i,\tau}-\tilde{t}_{i,\tau})(t_{i,\tau}+\tilde{t}_{i,\tau}-2\tau)-c^2(t_{i+1,\tau}-\tilde{t}_{i+1,\tau})(t_{i+1,\tau}+\tilde{t}_{i+1,\tau}-2\tau),
\end{align*}
for $i=1,3,5$. 
Furthermore, knowing that
$$r({\bf x_j},t_{j,\tau})=\tilde{r}({\bf x_j},\tilde{t}_{j,\tau})=\tau\ \forall 1\leq j\leq 6,$$
we get
\begin{align*}
{t}_{j,\tau}-\tilde{t}_{j,\tau}& = {t}_{j,\tau}-\tilde{r}^{-1}({\bf x_j},\tau)\\
& = {t}_{j,\tau}-\tilde{r}^{-1}({\bf x_j},r({\bf x_j},{t}_{j,\tau}))\\
& = \tilde{r}^{-1}({\bf x_j},\tilde{r}({\bf x_j},{t}_{j,\tau}))-\tilde{r}^{-1}({\bf x_j},r({\bf x_j},{t}_{j,\tau}))\\
&= \tilde{r}^{-1}\left({\bf x_j},\tilde{r}({\bf x_j},{t}_{j,\tau})-r({\bf x_j},{t}_{j,\tau})\right).
\end{align*}
Therefore,
\begin{align*}
|{t}_{j,\tau}-\tilde{t}_{j,\tau}|&\leq \|(\tilde{r}^{-1})'\|_{L^\infty(0,T)}|\tilde{r}({\bf x_j},{t}_{j,\tau})-r({\bf x_j},{t}_{j,\tau})|.
\end{align*}
Moreover, we deduce from \eqref{r} that $(\tilde{r}^{-1})'= \tilde{h}({\bf x_j},\tilde{r})$, which implies that
\begin{align*}
|{t}_{j,\tau}-\tilde{t}_{j,\tau}|&\leq \|\tilde{h}({\bf x_j},.)\|_{L^\infty(0,T)}|\tilde{r}({\bf x_j},{t}_{j,\tau})-r({\bf x_j},{t}_{j,\tau})|\\
&\leq 2 \|\tilde{r}({\bf x_j},.)-r({\bf x_j},.)\|_{L^\infty(0,T_{\textrm{obs}})}.
\end{align*}
Following the results of Theorem \ref{r estimate}, we finally deduce that 
$${t}_{j,\tau}-\tilde{t}_{j,\tau}\leq C\|\widetilde{\phi}({\bf x_j,.})-\phi({\bf x_j,.})\|_{L^\infty(0,T_{\textrm{obs}})}.$$
Therefore,
\begin{align*}
	A_{i\ i+1}-\widetilde{A}_{i\ i+1}\leq
	& 2c^2 (T_{\textrm{obs}}-\tau)C\\&\times \left(\|\widetilde{\phi}({\bf x_i,.})-\phi({\bf x_i,.})\|_{L^\infty(0,T_{\textrm{obs}})}+\|\widetilde{\phi}({\bf x_{i+1},.})-\phi({\bf x_{i+1},.})\|_{L^\infty(0,T_{\textrm{obs}})}\right).
\end{align*}
Finally, we deduce from \eqref{b} that $\forall \tau\in (0,T)$
\begin{align*}
|b(\tau)-\tilde{b}(\tau)|\leq& c^2 (T_{\textrm{obs}}-\tau)C\|X^{-1}\|_\infty\\&
\times\underset{i=1,3,5}{\sup}\left(\|\widetilde{\phi}({\bf x_i,.})-\phi({\bf x_i,.})\|_{L^\infty(0,T_{\textrm{obs}})}+\|\widetilde{\phi}({\bf x_{i+1},.})-\phi({\bf x_{i+1},.})\|_{L^\infty(0,T_{\textrm{obs}})}\right),
\end{align*}
which implies \eqref{b estimate}.
\end{preuve}

\section{Conclusion}
In this paper  we provide   a reconstruction procedure of the trajectory of a point source for the wave equation 
from the knowledge of the field at six well-chosen points on the observation boundary over a finite time interval.  
We derived a Lipschitz  stability estimate for the inversion that shows that the inverse problem is in fact well 
posed. The method can be  easily extended to other dimensions.  In this case the reconstruction requires 
measurement of the field at $2d$ well chosen points where $d$ is the dimension of the space. We plan in the
 future works  to numerically  implement the  method,  and  to recover simultaneously the intensity of the point
  source.  
\section*{}
The present work has been realized as part of Multi-Ondes projects receiving financial supports from ANR-France


\begin{thebibliography}{10}
\bibitem{Ammari02}   Ammari, H.,  Bao, G., \& Fleming, J. L. (2002). An inverse source problem for Maxwell's equations in magnetoencephalography. SIAM Journal on Applied Mathematics, 62(4), 1369-1382.
ISO 690	
%
\bibitem{Bao10} Bao, G., Lin, J., \& Triki, F. (2010). A multi-frequency inverse source problem. Journal of Differential Equations, 249(12), 3443-3465.
%
\bibitem{Bao11} Bao, G., Lin, J., \& Triki, F. (2011). Numerical solution of the inverse source problem for the Helmholtz equation with multiple frequency data. Contemp. Math, 548, 45-60.
%
\bibitem{Bao21}  Bao, G., Liu, Y., \& Triki, F. (2021). Recovering point sources for the inhomogeneous Helmholtz equation. Inverse Problems, 37(9), 095005.
%
\bibitem{Bao212} Bao, G., Liu, Y., \& Triki, F. Recovering Simultaneously a Potential and a Point Source from Cauchy Data.
Minimax Theory Appl.  6, No. 2, 227-238 (2021).
%
\bibitem{Bruncker2000} Bruckner, G., \& Yamamoto, M. (2000). Determination of point wave sources by pointwise observations: stability and reconstruction. Inverse problems, 16(3), 723.
%
\bibitem{Badia2001} El Badia, A., \& Ha-Duong, T. (2001). Determination of point wave sources by boundary measurements. Inverse Problems, 17(4), 1127.
%
\bibitem{Badia2002} El Badia, A., \& Ha-Duong, T. (2002). On an inverse source problem for the heat equation. Application to a pollution detection problem. Journal of inverse and ill-posed problems, 10(6), 585-599.
%
\bibitem{Hu}Hu, G., Kian, Y., Li, P., \& Zhao, Y. (2019). Inverse moving source problems in electrodynamics. Inverse Problems, 35(7), 075001.
%
\bibitem{Isakov90} Isakov, V. (1990). Inverse source problems (No. 34). American Mathematical Soc..
%
\bibitem{Jackson99} Jackson J.D. (1999). Classical Electrodynamics 3rd edn (NewYork:Wiley).
%
\bibitem{Komornik2002}Komornik V, \& Yamamoto M. Upper and lower estimates in determining point sources in a wave equation. Inverse Problems. 2002;18(2):319–329.
%
\bibitem{Komornik2005} Komornik, V., \& Yamamoto, M. (2005). Estimation of point sources and applications to inverse problems. Inverse Problems, 21(6), 2051.
%
\bibitem{Liu}Liu, Y., Hu, G., \& Yamamoto, M. (2021). Inverse moving source problem for time-fractional evolution equations: determination of profiles. Inverse Problems, 37(8), 084001.
%
\bibitem{Nakaguchi2012} Nakaguchi, E., Inui, H., \& Ohnaka, K. (2012). An algebraic reconstruction of a moving point source for a scalar wave equation. Inverse Problems, 28(6), 065018.
%
\bibitem{Ohe2011}Ohe, T., Inui, H., \& Ohnaka, K. (2011). Real-time reconstruction of time-varying point sources in a three-dimensional scalar wave equation. Inverse Problems, 27(11), 115011.
%
\bibitem{Ohe2020}Ohe, T. (2020). Real-time reconstruction of moving point/dipole wave sources from boundary measurements. Inverse Problems in Science and Engineering, 28(8), 1057-1102.
%
\bibitem{Rashedi2015}Rashedi, K., \& Sini, M. (2015). Stable recovery of the time-dependent source term from one measurement for the wave equation. Inverse Problems, 31(10), 105011.
\end{thebibliography}
\end{document}